\documentclass[11pt]{amsart}
\usepackage{amsfonts}
\usepackage{amsfonts}
\usepackage{amsfonts}
\usepackage{amsfonts}
\theoremstyle{plain}
\newtheorem{Thm}{Theorem}

\newtheorem{Lem}[Thm]{Lemma}

\newtheorem{Rk}[Thm]{Remark}

\begin{document}
\large

\title[blow-up]
{Blow-up of solutions to the nonlinear Schr$\bf{\ddot{O}}$dinger
equations on standard N-sphere and hyperbolic N-space}
\author{ Li Ma and Lin Zhao}
\address{Li Ma, Department of Mathematical Sciences, Tsinghua University,
 Peking 100084, P. R. China}

\email{lma@math.tsinghua.edu.cn}

\thanks{The research is partially supported by the National Natural Science
Foundation of China 10631020 and SRFDP 20060003002}

\begin{abstract}
In this paper, we partially settle down the long standing open
problem of the finite time blow-up property about the nonlinear
Schr$\ddot{o}$dinger equations on some Riemannian manifolds like the
standard 2-sphere $S^2$ and the hyperbolic 2-space $H^{2}(-1)$.
Using the similar idea, we establish such blow-up results on higher
dimensional standard sphere and hyperbolic $n$-space. Extensions to
$n$-dimensional Riemannian warped product manifolds with $n\geq 2$
are also given.


\end{abstract}

\maketitle
\date{26-6-2005}

\section{Introduction}
In this paper, we partially settle down the long standing open
problem of the finite time blow-up property about the nonlinear
Schr$\ddot{o}$dinger equations on some Riemannian manifolds like the
standard 2-sphere $S^2$ and the hyperbolic 2-space $H^{2}(-1)$. The
nonlinear Schr$\ddot{o}$dinger equations of the following form
\begin{equation}
iu_t=\Delta u+F(|u|^2)u
\end{equation}
play an important role in many areas of applied physics, such as
non-relativistic quantum mechanics, laser beam propagation,
Bose-Einstein condensates and so on (see \cite{Su}). The initial
value problems (IVP) or the initial-boundary value problems (IBVP)
of (1) on $\mathbb{R}^n$ have been extensively studied in the last
two decades (see \cite{CW, Kat, Y, GV, We1, We2} ). In particular,
the blow-up properties in finite time for IVP or IBVP have caught
sufficient attention (see \cite{Gl, Ka, Og1, Me1, Me2}). However,
much less results have been known on bounded domains in
$\mathbb{R}^n$ or on compact manifolds $(M,g)$, with the notable
exception of the works of H.Br$\acute{e}$zis and T.Gallouet
\cite{Br}, J.Bourgain \cite{Bo1, Bo2, Bo3} (In \cite{Bo1}, the case
$M=\mathbb{R}^2/\mathbb{Z}^2$ was discussed in detail), and N.Burq,
P.G$\acute{e}$rard and N.Tzvetkov \cite{BGT1, BGT2}. In particular,
the blow-up in finite time of Schr$\ddot{o}$dinger equations (1)
posed on an arbitrary Riemannian manifold ($M$,$g$) is a long widely
known open problem. To our knowledge, the only examples of such
blow-up phenomena on Riemannian manifolds are given by the following
result, attributed to Ogawa-Tsutsuni \cite{Og2} if the dimension $n$
of $M$ equals to 1 and generalized to the case $n=2$ by N.Burq,
P.G$\acute{e}$rard and N.Tzvetkov \cite{BGT1}.

\begin{Thm}
Let $(M,g)$ be a compact Riemannian manifold of dimension $n=1$ or
$n=2$. Assume there exist $x^0\in M$ and a system of coordinates
near $x^0$ in which
$$
g=\sum_{j=1}^d dx_j^2.
$$
Then there exist smooth solutions $u\in C^\infty([0,T)\times M)$ of
$$
iu_t=\Delta u+|u|^{4/n}u
$$
such that as $t\rightarrow T$,
$$
|u(t,x)|^2\rightharpoonup
\|Q\|^2_{L^2({\mathbb{R}}^n)}\delta(x-x^0),
$$
where $Q$ is the ground state solution on $\mathbb{R}^n$ of
$$
\Delta Q+Q^{1+4/n}=Q.
$$
\end{Thm}

Even though the condition of Theorem 1 that the manifold near $x^0$
is flat is a very strong restriction, the result is also impressive.

In this paper, we concentrate on the analysis of the blow-up
phenomena for IVP or BIVP of the Schr$\ddot{o}$dinger equations
posed on Riemannian manifolds. To be precise, the IVP and BIVP are
of the following forms respectively
\begin{align}
\textrm{IVP}\left\{\begin{array}{ll}
iu_t=\Delta u+F(|u|^2)u,\ \textrm{on}\ M, \\
u(0,x)=u_0(x),\\
\partial M=\emptyset;
\end{array}
\right.
\end{align}
\begin{align}
\textrm{BIVP}\left\{\begin{array}{ll}
iu_t=\Delta u+F(|u|^2)u,\ \textrm{on}\ M, \\
u|_{{\bf{R}}\times\partial M}=0,\\
u(0,x)=u_0(x),
\end{array}
\right.
\end{align}
where $F$ is a real-valued smooth function on the n-dimensional
Riemannian manifold $(M,g)$ and $F$ satisfies $F(s)\leq
C(1+s^{(p-1)/2})$ for some $p>1$ on $[0,\infty)$. Here $\Delta$ is
the Laplacian operator of the metric $g$ with the sign $\Delta
u=u^{''}$ on the real line $\mathbb{R}$. Noticing that (3) reduces
to (2) when $\partial M=\emptyset$, it's convenient to establish our
blow-up results in the context of (3). The local wellposedness in
$H^s(M)$ for $s>n/2$ to (3) with the interesting case
$F(|u|^2)u=|u|^{p-1}u$ ($p>1$) is a classical consequence of energy
estimates, and therefore, it's relaxed to assume that the solution
$u(t)$ to (3) satisfies
\begin{align}
u\in C^1([0,T),L^2)\cap C([0,T), H^2\cap H^1_0\cap L^{p+1}),
\end{align}
where $T$ is the maximal existence time for the solution $u(t)$.

Before stating our main results, let's introduce some exact
notations concerning Riemannian manifolds used below. Let $(M,g)$ be
a complete Riemannian manifold of dimension $n$ with boundary
$\partial M$ or not. We denote by $D$ the Levi-Civita connection,
and by $TM=\bigcup_{x\in M}T_xM$, where $T_xM$ is the tangent space
at $x\in M$. It's well known that the smooth sections of $TM$ are
just vector fields. For $f\in C^1(M)$, its gradient is defined as
the unique vector field $\nabla f$ such that
$$
\forall x\in M,\ \forall \xi\in TM, \ \ g(\nabla f(x),\xi(x))=(\xi
f)(x).
$$
The divergence $\textrm{div} \bf{X}$ of a smooth vector field
$\bf{X}$ is defined as the unique smooth function on $M$ such that
$$
\forall f\in C^\infty_0(M),\ \ \int f\textrm{div}{\bf{X}}=-\int
{\bf{X}}f.
$$
The Laplace-Beltrami operator $\Delta$ on $M$ is the second order
differential operator defined by
$$
\forall f\in C^2(M),\ \ \Delta f=\textrm{div}(\nabla f).
$$
Corresponding to our analysis, we need to extend $g$ to be defined
on complex valued vector fields. For complex valued vectors
$X_1+iX_2$, $Y_1+iY_2$, where $X_1$, $X_2$, $Y_1$, $Y_2$ are real,
we define
\begin{align*}
&\quad g(X_1+iX_2,Y_1+iY_2)\\
&=g(X_1,X_2)+ig(X_1,Y_1)+ig(X_2,Y_1)+g(Y_1,Y_2).
\end{align*}
It's easy to see that $g$ defined in such a way is bilinear in the
field $\mathbb{C}$ and accordingly
$$
\nabla f=\nabla\Re f+i\nabla\Im f,\ \ \Delta f=\Delta\Re
f+i\Delta\Im f.
$$
When $M$ has a nonempty boundary $\partial M$, we denote by
${\bf{v}}$ the outer unit normal vector along $\partial M$.

For the sake of simplicity, we omit the spatial integral variable
$x\in M$ and omit the integral region when it's the whole space $M$,
and we abbreviate $L^q(M)$, $H^k(M)$ to $L^q$, $H^k$ respectively.
We write the integral $\int_M dV_M$ and $\int_{\partial
M}dV_{\partial M}$ as $\int$ and $\oint$ respectively, and the norm
of $L^q$ as $\|\cdot\|_q$. We denote by $S^n$ the standard sphere of
dimension $n$ and by $H^n(-1)$ the hyperbolic $n$-space
respectively. We denote by $\textrm{N}$ the north pole of $S^n$ and
by $\textrm{dist}(\textrm{N},x)$ the distance between $\textrm{N}$
and $x\in S^n\setminus\{\textrm{N}\}$.

Our main results in 2-dimensions are the following two Theorems.

\begin{Thm}
Consider the Schr$\ddot{o}$dinger equation
\begin{align*}
\left\{\begin{array}{ll}
iu_t=\Delta u+|u|^{p-1}u,\ \ \textrm{on}\ S^2,\\
u(0)=u_0\in H^{1}.
\end{array}
\right.
\end{align*}
For $p\geq5$, if $u_0(x)=u_0(r)$, where
$r=\textrm{dist}(\textrm{N},x)$ and $u_0(r)$ is an asymmetric
function at $r=\frac{\pi}{2}$ with $u(\frac{\pi}{2})=0$ and $E_0<0$,
then the asymmetric solution satisfying (4) blows up in finite time.
\end{Thm}

\begin{Thm}
Consider the Schr$\ddot{o}$dinger equations (2) on $M=H^2(-1)$.
Assume there exists a constant $\kappa\geq3$ such that
$$
sF(s)\geq \kappa G(s),\ \ \forall s\geq0.
$$
Then any solution satisfying (4) with $E_0<0$ blows up in finite
time.
\end{Thm}

For higher dimensional results, please see Theorem 12 and 13. Here,
we just want to point out that the range of the exponent $p$ for
blow-up of the solutions to the Schr$\ddot{o}$dinger equation (1)
with $F(s)=s^{(p-1)/2}$ on $\mathbb{R}^n$ is $p\geq 1+\frac{4}{n}$,
on $S^n$ is $p\geq 5$ and on $H^{n}(-1)$ is $p\geq 1+\frac{4}{n-1}$.

We try to present very elementary proofs of our blow up results
starting from 2-sphere. This paper is organized as follows. In
section 2, we establish some new invariant quantities for the
Schr$\ddot{o}$dinger equations on general Riemannian manifolds,
which generalize the corresponding classical results on
$\mathbb{R}^n$. In section 3, we construct blow-up solutions on the
unit sphere $S^2$. In section 4, we establish the blow-up results on
a class of noncompact manifolds. We discuss the blow-up properties
on $n$-dim manifolds with $n\geq3$ in section 5.

\section{Preliminary lemma}
The following lemma is a generalization of the identities obtained
by Glassey \cite{Gl} (see also \cite{Ka}). We define
$$
G(u)=\int_0^uF(s)ds.
$$
\begin{Lem}
Suppose that $(M,g)$ is a complete Riemannian manifold of dimension
$n$ with boundary $\partial M$ or not, and ${\bf{v}}$ is the outer
unit normal vector along $\partial M$. Let $u$ be a solution of (3)
satisfying (4), $\rho$ be an arbitrary smooth function on $M$,
and $\bf{X}$ be a real smooth vector field on $M$. Define $J(t):=\int \rho|u|^2$. Then we have\\
(A). $\|u(t)\|_2=\|u_0\|_2$,\\
(B). $\int(g(\nabla u,\nabla\bar{u})-G(|u|^2))\equiv \textrm{const}:=E_0$,\\
(C). $J'(t)=-2\Im\int g(\nabla\rho,\nabla u)\bar{u}$,\\
(D).
\begin{align*}
&\quad\frac{d}{dt}\Im\int g({\bf{X}},\nabla u)\bar{u}\\
&=-2\int
D{\bf{X}}(\nabla u,\nabla\bar{u})+
\frac{1}{2}\int(\Delta \textrm{div}{\bf{X}})|u|^2\\
&\quad+\int (\textrm{div}{\bf{X}})(F(|u|^2)|u|^2-G(|u|^2))\\
&\quad+\oint g(\nabla u,\nabla \bar{u})g({\bf{X}},{\bf{v}}).
\end{align*}
\end{Lem}

\begin{proof}
The facts that $\int \rho|u|^2$ and $\Im\int g({\bf{X}},\nabla
u)\bar{u}$ are of $C^1[0,T)$ are straightforward and the reader can
refer to \cite{Ka} for details.

For (A), multiply both sides of (3) by $2\bar{u}$ and take the
imaginary part to obtain
\begin{align}
\frac{\partial}{\partial t}|u|^2=2\nabla\cdot \Im(\bar{u}\nabla u).
\end{align}
Integrating it over $M$ we get (A).

For (B), multiply (3) by $2\bar{u}_t$, integrate, and take the real
part of the resulting expression.

For (C), multiply (5) by $\rho$ and integrate by parts over $M$.

The derivation of (D) is a bit involved. We first multiply (3) by
$2D_{\bf{X}}\bar{u}$ to obtain
\begin{align}
2i(D_{\bf{X}}\bar{u})u_t&=2(D_{\bf{X}}\bar{u})\Delta
u+2(D_{\bf{X}}\bar{u})F(|u|^2)u\\
&:=\textrm{I}_1+\textrm{I}_2.\nonumber
\end{align}
Then, we take the real part of the left-hand side (LHS) of (6) to
get
\begin{align*}
\Re
(\textrm{LHS})&=i((D_{\bf{X}}\bar{u})u_t-(D_{\bf{X}}u)\bar{u}_t)\\
&=i((uD_{\bf{X}}\bar{u})_t-D_{\bf{X}}(u\bar{u}_t))\\
&=\Re(i(uD_{X}\bar{u})_t)-\Re(iD_{\bf{X}}(u\bar{u}_t))\\
&=\frac{d}{dt}\Im (g({\bf{X}},\nabla u)\bar{u})-\Re
(iD_{\bf{X}}(u\bar{u}_t)).
\end{align*}
Integrating this identity over $M$ yields
\begin{align}
\Re\int \textrm{LHS} =\frac{d}{dt}\Im\int g({\bf{X}},\nabla
u)\bar{u}-\Re\int iD_{\bf{X}}(u\bar{u}_t).
\end{align}
Using integration by parts we have
\begin{align*}
&\quad\Re \int iD_{\bf{X}}(u\bar{u}_t)\\
&=-\Re\int
(\textrm{div}{\bf{X}})(iu\bar{u}_t)\\
&=-\Re\int (\textrm{div}{\bf{X}})(-u\Delta\bar{u}-F(|u|^2)|u|^2)\\
&=-\Re\int(g(\nabla
(\textrm{div}{\bf{X}}),u\nabla\bar{u})+(\textrm{div}{\bf{X}})g(\nabla
u,\nabla \bar{u}))+
\int (\textrm{div}{\bf{X}})F(|u|^2)|u|^2\\
&=-\frac{1}{2}\int
g(\nabla(\textrm{div}{\bf{X}}),\nabla|u|^2)-\int(\textrm{div}{\bf{X}})g(\nabla
u,\nabla\bar{u})\\
&\quad+\int (\textrm{div}{\bf{X}})F(|u|^2)|u|^2\\
&=\frac{1}{2}\int(\Delta\textrm{div}{\bf{X}})|u|^2-\int(\textrm{div}{\bf{X}})g(\nabla
u,\nabla\bar{u})+\int (\textrm{div}{\bf{X}})F(|u|^2)|u|^2.
\end{align*}
Inserting this into (7) we obtain that
\begin{align}
\Re\int \textrm{LHS} &=\frac{d}{dt}\Im\int g({\bf{X}},\nabla
u)\bar{u}-\frac{1}{2}\int(\Delta\textrm{div}{\bf{X}})|u|^2\\
&\quad+\int(\textrm{div}{\bf{X}})g(\nabla u,\nabla
\bar{u})-\int(\textrm{div}{\bf{X}})F(|u|^2)|u|^2.\nonumber
\end{align}

To handle the right-hand side of (6), we use integration by parts
again to obtain
\begin{align*}
\Re\int \textrm{I}_1&=2\Re\int (D_{\bf{X}}\bar{u})\Delta u\\
&=-2\Re\int g(\nabla D_{\bf{X}}\bar{u},\nabla u)+2\Re\oint
g((D_{\bf{X}}\bar{u})\nabla u,{\bf{v}})\\
&=-2\Re\int (D{\bf{X}}(\nabla u,\nabla
\bar{u})+\frac{1}{2}D_{\bf{X}}g(\nabla u,\nabla\bar{u}))\\
&\quad+2\Re\oint g((D_{\bf{X}}\bar{u})\nabla u,{\bf{v}})\\
&=-2\int D{\bf{X}}(\nabla u,\nabla \bar{u})-\int
D_{\bf{X}}g(\nabla u,\nabla\bar{u})\\
&\quad+2\oint g(\nabla u,\nabla\bar{u})g({\bf{X}},{\bf{v}}).
\end{align*}
The last "=" in the above expression follows from the fact that
$u|_{\partial M}=0$, which implies $\nabla u=g(\nabla
u,{\bf{v}}){\bf{v}}$.

Noticing that
\begin{align*}
\int D_{\bf{X}}g(\nabla u,\nabla\bar{u})&=-\int
(\textrm{div}{\bf{X}})g(\nabla u,\nabla\bar{u})+\oint g(\nabla
u,\nabla\bar{u})g({\bf{X}},{\bf{v}}),
\end{align*}
we then get
\begin{align}
\Re\int \textrm{I}_1&=-2\int D{\bf{X}}(\nabla u,\nabla \bar{u})+\int
(\textrm{div}{\bf{X}})g(\nabla u,\nabla\bar{u})\\
&\quad+\oint g(\nabla u,\nabla\bar{u})g({\bf{X}},{\bf{v}}).\nonumber
\end{align}
For $I_2$ we have
\begin{align}
\Re\int \textrm{I}_2&=2\Re\int (D_{\bf{X}}\bar{u})F(|u|^2)u=\int
(D_{\bf{X}}|u|^2)F(|u|^2)\\
&=\int D_{\bf{X}}G(|u|^2)=-\int
(\textrm{div}{\bf{X}})G(|u|^2).\nonumber
\end{align}

Combining (8)-(10) with (6) we get (D), and the proof of the lemma
is concluded.
\end{proof}

\begin{Rk}
If we choose ${\bf{X}}=\nabla \rho$ in lemma 2 (D), we then arrive
at
\begin{align}
J^{''}(t)&=4\int D^2\rho(\nabla u,\nabla\bar{u})-\int
(\Delta^2\rho)|u|^2\\
&\quad -2\int(\Delta\rho)(F(|u|^2)|u|^2-G(|u|^2))\nonumber\\
&\quad -2\oint g(\nabla
u,\nabla\bar{u})g(\nabla\rho,{\bf{v}}).\nonumber
\end{align}
This identity will play a vital role in our analysis. In particular,
when $\nabla\rho={\bf{v}}$, we have
$$
g(\nabla u,\nabla\bar{u})g(\nabla\rho,{\bf{v}})=g(\nabla
u,\nabla\bar{u})\geq0,
$$
which is an important fact in our proof.
\end{Rk}

\section{Blow-up phenomena on $S^2$}
To investigate the blow-up nature of the solution $u$, one method is
to observe the long time behavior of $J(t):=\int\rho|u|^2$. If there
exists $\rho\geq0$ such that $J(t)$ becomes negative after some
finite time $T$ due to the conservation of the $L^2$ norm and the
energy $E_0$, then $u$ must blow up before the time $T$. It's
classical that $\rho(x)=|x|^2$ when $M=\mathbb{R}^d$. But for an
arbitrary manifold, the sharp $\rho$ adopted to the blow-up
properties is unknown explicitly. It seems that $\Delta |x|^2=2n$ is
a nice property for us to use (11) on $M=\mathbb{R}^n$. For
noncompact manifolds and compact manifolds with boundary, it's
possible to find $\rho$ such that $\Delta\rho=\textrm{const}$. We
put this idea in practice on the half sphere
$S^2_{+}:=S^2\cap\{x^1\geq0\}$ and the hyperbolic space $M=H^2(-1)$.
We now state the result for $S^2_{+}$.

\begin{Thm}
Consider the Schr$\ddot{o}$dinger equations (3) on $M=S^2_{+}$.
Assume there exists a constant $\kappa\geq3$ such that
$$
sF(s)\geq \kappa G(s),\ \ \forall s\geq0.
$$
Then any solution satisfying (4) with $E_0<0$ blows up in finite
time.
\end{Thm}

\begin{proof}
In $\mathbb{R}^3$,
$S^2_{+}=\{(x^1)^2+(x^2)^2+(x^3)^2=1\}\cap\{x^1\geq0\}$. We want to
construct a function $\rho$ such that

(a). $\rho\in C^4(S^2_{+})$;

(b). $\rho>0$ on $S^2_{+}\setminus\{\textrm{N}\}$;

(c). $\Delta \rho=1$ on $S^2_{+}$;

(d). $D^2\rho(\nabla u,\nabla\bar{u})\leq g(\nabla
u,\nabla\bar{u}),\ \forall u\in C^1(S^2_{+})$.

We now give the form of $\rho$ exactly. To make the calculations
clear, we recall the expressions of $\nabla f$, $\Delta f$ and
$D^2f$ for $f\in C^2(M)$ in local coordinates. We use Einstein's
convention. Let $g=g_{ij}dx^idx^j$, $G=\det(g_{ij})$ and
$(g^{ij})=(g_{ij})^{-1}$. Then
\begin{align}
\nabla f=g^{ij}f_i\partial_j,
\end{align}
\begin{align}
\Delta f=\frac{1}{\sqrt{G}}\partial_i(g^{ij}\sqrt{G}f_j),
\end{align}
\begin{align}
D^2f=(f_{ij}-\Gamma_{ij}^kf_k)dx^i\otimes dx^j,
\end{align}
where
$$
\Gamma_{ij}^k=\frac{1}{2}g^{kl}(\frac{\partial g_{il}}{\partial
x^j}+\frac{\partial g_{lj}}{\partial x^i}-\frac{\partial
g_{ij}}{\partial x^l}).
$$

In our situation, we use the geodesic polar coordinate $(r,\theta)$
at the north pole $\textrm{N}$ for $S^2_{+}$, i.e.,
\begin{align*}
\left\{\begin{array}{ll}
x^1=\cos r,\\
x^2=\sin r\cos\theta,\\
x^3=\sin r\sin\theta,
\end{array}
\right.
\end{align*}
where $r\in(0,\pi/2]$, $\theta\in[0,2\pi)$. In this coordinate,
\begin{align*}
\left\{\begin{array}{ll}
dx^1=-\sin rdr,\\
dx^2=\cos r\cos\theta dr-\sin r\sin\theta d\theta,\\
dx^3=\cos r\sin\theta dr+\sin r\cos\theta d\theta,
\end{array}
\right.
\end{align*}
and hence
$$
g=\sum_{i=1}^3 (dx^i)^2=dr^2+\sin^2rd\theta^2,
$$
$$
\Gamma_{11}^1=\Gamma_{12}^1=\Gamma_{21}^1=0,\ \ \Gamma_{22}^1=-\sin
r\cos r.
$$

By (13) we have for $\varrho=\varrho(r)\in C^2(0,\pi/2]$ that
\begin{align*}
\Delta \varrho(r)&=\frac{1}{\sin r}(\sin
r\varrho'(r))'=\varrho^{''}(r)+\varrho'(r)\cot r.
\end{align*}
Solving the ODE
\begin{align*}
\left\{\begin{array}{ll}
\varrho^{''}(r)+\varrho'(r)\cot r=1,\ \ 0<r\leq\pi/2,\\
\varrho(r)>0,\ \ 0<r\leq\pi/2,\\
\varrho'(\pi/2)=1,
\end{array}
\right.
\end{align*}
we get a solution $\varrho(r)=-2\log\cos(r/2)$. We then define
$\rho(r)$ as
\begin{align*}
\rho(r)=\left\{\begin{array}{ll}
\varrho(r),\ \ 0<r\leq\pi/2,\\
0,\ \ r=0.
\end{array}
\right.
\end{align*}
It's then easy to see that $\rho(r)\in C^4[0,\pi/2]$, i.e, $\rho\in
C^4(S^2_{+})$.

From (12) and (14), we have for any $u\in C^1(S^2_{+})$,
$$
g(\nabla u,\nabla\bar{u})=|u_{r}|^2+\frac{1}{\sin^2r}|u_{\theta}|^2
$$
and
\begin{align*}
D^2\rho(\nabla u,\nabla\bar{u})&=\rho^{''}(r)|u_{r}|^2
+\rho'(r)\frac{\cos r}{\sin^3r}|u_{\theta}|^2\\
&=\frac{1-\cos r}{\sin^2 r}|u_{r}|^2+
\frac{(1-\cos r)\cos r}{\sin^4r}|u_{\theta}|^2.\\
\end{align*}
It's obvious that $D^2\rho(\nabla u,\nabla\bar{u})\leq g(\nabla
u,\nabla\bar{u})$ provided  $0<r\leq\pi/2$. By a standard
approximation process we get
$$
\int_{S^2_{+}}D^2\rho(\nabla
u,\nabla\bar{u})\leq\int_{S^2_{+}}g(\nabla u,\nabla\bar{u}), \ \
\forall\ u\in H^1(S^2_{+}).
$$

Notice that $\nabla \rho=\bf{v}$ on the boundary $\partial S^2_{+}$.
Then by (11) we have
\begin{align}
J^{''}(t)&\leq 4\int_{S^2_{+}} D^2\rho(\nabla u,\nabla\bar{u})-2\int_{S^2_{+}}(F(|u|^2)|u|^2-G(|u|^2))\\
&\quad-2\oint_{\partial S^2_{+}}g(\nabla u,\nabla\bar{u})g(\nabla\rho,{\bf{v}})\nonumber\\
&\leq 4\int_{S^2_{+}} g(\nabla
u,\nabla\bar{u})-2(\kappa-1)\int_{S^2_{+}} G(|u|^2).\nonumber
\end{align}
Combining Lemma 2 (B) with (15) we obtain
$$
J^{''}(t)\leq 4E_0+(6-2\kappa)\int_{S^2_{+}} G(|u|^2)\leq 4E_0<0,
$$
which implies that $J(t)$ becomes negative after some finite time
$T$, i.e., the solution $u$ satisfying $E_0<0$ blows up in finite
time.
\end{proof}

\begin{Rk}
Here we only established the blow-up result for the
Schr$\ddot{o}$dinger equations on $S^2_{+}$. The method above can't
be applied to the whole sphere $S^2$ or some other compact manifolds
because of the nonexistence of subharmonic functions on compact
manifolds. But this result allows us to construct blow-up solutions
for the Schr$\ddot{o}$dinger equations posed on $S^2$. See the proof
of Theorem 2.
\end{Rk}

\begin{Rk}
If $F(s)=s^\frac{p-1}{2}$ where $p>1$, then
$G(s)=\frac{2}{p+1}s^{\frac{p+1}{2}}$, and the condition $sF(s)\geq
\kappa G(s)$ on $[0,+\infty)$ for some $\kappa\geq 3$ is equivalent
to $p\geq 5$. It's already known that when $1<p<3$, the
Schr$\ddot{o}$dinger equation $iu_t=\Delta u+|u|^{p-1}u$ with
$u_0\in H^s$ ($s\geq1$) on $S^2$ has a unique global solution $u\in
C(\mathbb{R},H^s)$ (see \cite{BGT2}). We point out when $p\geq5$,
blow-up phenomena may occur by our Theorem 2.
\end{Rk}

\emph{Proof of Theorem 2:} If $u_0(r)$ is asymmetric in $r$, then by
the symmetry of the Schr$\ddot{o}$dinger equation, the solution
$u(t,x)$ is also asymmetric, that is, $u(t,r)$ is asymmetric with
respected to $r$, which implies that $u(t,\pi/2)\equiv0$. We then
cut the sphere $S^2$ into two parts $S^2_{+}$ and $S^2_{-}$, where
$$
S^2_{+}=\{(x^1)^2+(x^2)^2+(x^3)^2=1\}\cap\{x^1\geq0\},
$$
$$
S^2_{-}=\{(x^1)^2+(x^2)^2+(x^3)^2=1\}\cap\{x^1\leq0\}.
$$
We write by $\rho_{+}$ the function obtained in the proof of Theorem
6 on $S^2_{+}$, and define
$$
\rho_{-}(x):=\rho_{+}(-x)\ \ \textrm{for} \ x\in S^2_{-}.
$$
We denote by ${\bf{v}}_{+}$ and ${\bf{v}}_{-}$ the outer normal
vector along $\partial S^2_{+}$ and $\partial S^2_{-}$. Then
${\bf{v}}_{+}=-{\bf{v}}_{-}$, and $\nabla\rho_{+}={\bf{v}}_{+}$,
$\nabla\rho_{-}={\bf{v}}_{-}$.

We now define
$$
J(t):=\int_{S^2_{+}}\rho_{+}|u|^2+\int_{S^2_{-}}\rho_{-}|u|^2:=J_1+J_2.
$$

Since $u=0$ on $\partial S^2_{+}=\partial S^2_{-}$, we can use (11)
directly for both $J_1$ and $J_2$ to get that
\begin{align}
J^{''}_1(t)&\leq 4\int_{S^2_{+}} D^2\rho_{+}(\nabla u,\nabla\bar{u})-2\int_{S^2_{+}}(F(|u|^2)|u|^2-G(|u|^2))\\
&\quad-2\oint_{\partial S^2_{+}}g(\nabla u,\nabla\bar{u})g(\nabla\rho_{+},{\bf{v}}_{+})\nonumber\\
&\leq 4\int_{S^2_{+}} g(\nabla
u,\nabla\bar{u})-2(\kappa-1)\int_{S^2_{+}} G(|u|^2),\nonumber
\end{align}
and
\begin{align}
J^{''}_2(t)&\leq 4\int_{S^2_{-}} D^2\rho_{-}(\nabla u,\nabla\bar{u})-2\int_{S^2_{-}}(F(|u|^2)|u|^2-G(|u|^2))\\
&\quad-2\oint_{\partial S^2_{-}}g(\nabla u,\nabla\bar{u})g(\nabla\rho_{-},{\bf{v}}_{-})\nonumber\\
&\leq 4\int_{S^2_{-}} g(\nabla
u,\nabla\bar{u})-2(\kappa-1)\int_{S^2_{-}} G(|u|^2).\nonumber
\end{align}
We obtain from (16)+(17) that
\begin{align*}
J^{''}(t)&\leq 4\int_{S^2} g(\nabla
u,\nabla\bar{u})-2(\kappa-1)\int_{S^2} G(|u|^2)\\
&=4E_0+(6-2\kappa)\int_{S^2}G(|u|^2)\\
&\leq 4E_0<0
\end{align*}
provided $\kappa\geq3$, i.e., $p\geq5$, which implies that the
solution $u$ blows up in finite time. This is the end of proof.

\begin{Rk}
When $3\leq p<5$, the blow-up property of the solution on $S^2$
leaves open.
\end{Rk}

\section{Blow-up on noncompact 2-dim manifolds}
By the method of introducing some proper weight function $\rho$, we
can get similar blow-up results for noncompact manifolds. We first
prove the result for the hyperbolic 2-space $H^{2}(-1)$, i.e.,
Theorem 3, and then generalize it to a class of noncompact
manifolds.

\emph{Proof of Theorem 3:} For $H^2(-1)$, it's standard that for
$s\in[0,\infty)$ and $\theta\in[0,2\pi)$,
$$
g=\frac{1}{1+s^2}ds^2+s^2d\theta^2.
$$
Choose $r(s)=\sinh^{-1}(s)$, then we have
$dr=\frac{1}{\sqrt{1+s^2}}ds$, and thus
$$
g=dr^2+\sinh^2(r)d\theta^2.
$$
Calculating by (13), we get for $\rho=\rho(r)\in C^2(0,\infty)$,
$$
\Delta \rho(r)=\frac{1}{\sinh r}(\sinh r
\rho'(r))'=\rho^{''}(r)+\rho'(r)\coth r.
$$
Solving the ODE
\begin{align*}
\left\{\begin{array}{ll}
\rho^{''}(r)+\rho'(r)\coth r=1,\ \ 0<r<\infty,\\
\rho(r)>0,\ \ 0<r<\infty,
\end{array}
\right.
\end{align*}
we get a solution $\rho(r)=2\log\cosh(\frac{r}{2})$. It's easy to
see that by defining $\rho(0)=0$, $\rho\in C^4[0,\infty)$, i.e.,
$\rho\in C^4(H^2(-1))$.

From (12) and (14) we have for any $u\in C^1(H^2(-1))$,
\begin{align*}
g(\nabla u,\nabla\bar{u})=|u_r|^2+\frac{1}{\sinh^2r}|u_\theta|^2,
\end{align*}
and
\begin{align*}
D^2\rho(\nabla
u,\nabla\bar{u})&=\rho^{''}(r)|u_r|^2+\frac{\cosh r}{\sinh^3r}|u_\theta|^2\\
&=\frac{\cosh r-1}{\sinh^2r}|u_r|^2+\frac{\cosh r(\cosh r-1)}
{\sinh^4r}|u_\theta|^2.
\end{align*}

Noticing that when $r>0$ we have
$$
\frac{\cosh r-1}{\sinh^2r}=\frac{1}{2\cosh^2r/2}\leq\frac{1}{2}
$$
and
$$
\frac{\cosh r(\cosh r-1)}{\sinh^2r}=\frac{\cosh
r}{2\cosh^2r/2}\leq1,
$$
we obtain that
$$
D^2\rho(\nabla u,\nabla\bar{u})\leq g(\nabla u,\nabla\bar{u}).
$$
The remainder of the proof is the same as in Theorem 6. This
completes the proof of Theorem 3.

The next result generalize our analysis to a class of 2-dim
noncompact manifolds.

\begin{Thm}
Let $(M,g)$ be a 2-dim Riemannian manifold such that $M$ can be
covered by only one coordinate system in which
$$
g=dr^2+h^2(r)d\theta^2,
$$
where $h\in C^4[0,\infty)$ satisfying for some constants
$\tau_1,\tau_2\in[0,1]$,
\begin{align*}
\left\{\begin{array}{ll}
h(r)\in C^4[0,\infty),\\
h(r)>0,\ \textrm{on}\ (0,\infty); \ \ h(0)=0; \ \ h'(0)>0,\\
\tau_1h^2(r)\leq h'(r)\int_0^rh(s)ds\leq\tau_2h^2(r),\ \textrm{on}\
(0,\infty).
\end{array}
\right.
\end{align*}
Consider the Schr$\ddot{o}$dinger equations (2) on $M$. Assume there
exists a constant $\kappa\geq2\max\{1-\tau_1,\tau_2\}+1$ such that
$$
sF(s)\geq \kappa G(s),\ \ \forall s\geq0.
$$
Then any solution satisfying (4) with $E_0<0$ blows up in finite
time.
\end{Thm}

\begin{proof}
In above coordinate,
$$
\Gamma_{11}^1=\Gamma_{12}^1=\Gamma_{21}^1=0,\ \
\Gamma_{22}^1=-h(r)h'(r).
$$

By (13) we have for $\rho=\rho(r)\in C^2(0,\infty)$ that
\begin{align*}
\Delta
\rho(r)&=\frac{1}{h(r)}(h(r)\rho'(r))'=\rho^{''}(r)+\rho'(r)\frac{h'(r)}{h(r)}.
\end{align*}
Solving the ODE
\begin{align*}
\left\{\begin{array}{ll}
\rho^{''}(r)+\rho'(r)\frac{h'(r)}{h(r)}=1,\ \ 0<r<\infty,\\
\rho(r)>0,\ \ 0<r<\infty,
\end{array}
\right.
\end{align*}
we get a solution
$$
\rho(r)=\int_0^r(\int_0^sh(t)dt)(h(s))^{-1}ds.
$$
We then define $\rho(0)=0$, and thus it's easy to see that
$\rho(r)\in C^4[0,\infty)$, i.e, $\rho\in C^4(M)$.

From (12) and (14), we have for any $u\in C^1(M)$,
$$
g(\nabla u,\nabla\bar{u})=|u_{r}|^2+\frac{1}{h^2(r)}|u_{\theta}|^2
$$
and
\begin{align*}
D^2\rho(\nabla u,\nabla\bar{u})&=\rho^{''}(r)|u_{r}|^2
+\rho'(r)\frac{h'(r)}{h^3(r)}|u_{\theta}|^2\\
&=(1-\frac{h'(r)\int_0^rh(s)ds}{h^2(r)})|u_{r}|^2+
\frac{h'(r)\int_0^rh(s)ds}{h^4(r)}|u_{\theta}|^2\\
&\leq(1-\tau_1)|u_r|^2+\frac{\tau_2}{h^2(r)}|u_\theta|^2\\
&\leq\max\{1-\tau_1,\tau_2\}g(\nabla u,\nabla\bar{u}).
\end{align*}
The remainder of the proof is the same as in Theorem 6.
\end{proof}

\section{Some results for higher dimensional manifolds}
In this section, we make some calculations for $S^n$ and $H^n(-1)$.
We first compute the case $S^n_{+}=S^2\cap\{x^1\geq0\}$. As above,
for $0<r\leq\pi/2$, $0\leq\theta_1,\cdots,\theta_{n-1}<2\pi$, we
have
\begin{align*}
\left\{\begin{array}{ll}
x^1=\cos r,\\
x^2=\sin r\cos\theta_1,\\
x^3=\sin r\sin\theta_1\cos\theta_2,\\
\cdots\ \ \cdots\\
x^n=\sin r\sin\theta_1\cdots\sin\theta_{n-2}\cos\theta_{n-1},\\
x^{n+1}=\sin r\sin\theta_1\cdots\sin\theta_{n-2}\sin\theta_{n-1}.
\end{array}
\right.
\end{align*}
By (13) we get
\begin{align*}
g=dr^2+\sin^2r(d\theta_1^2+\sin^2\theta_1d\theta_2^2+\cdots+\sin^2\theta_1\cdots\sin^2\theta_{n-2}d\theta_{n-1}^2),
\end{align*}
and
\begin{align*}
\Delta \rho(r)=\frac{(\rho'(r)\sin^{n-1}r)'}{\sin^{n-1}r}.
\end{align*}
Solving the ODE
$$\Delta \rho(r)=1\ \ on \ \ (0,\pi/2]$$
we get the desired positive solution
\begin{align*}
\rho(r)=\int_0^r(\int_0^s\sin^{n-1}\tau d\tau)(\sin^{n-1}s)^{-1}ds.
\end{align*}
By defining $\rho(0)=0$ we see that $\rho\in C^4(S^n_{+})$.

From (12) and (14) we get for all $u\in C^1(S^n_{+})$
\begin{align*}
g(\nabla
u,\nabla\bar{u})&=|u_r|^2+\frac{|u_{\theta_1}|^2}{\sin^2r}+\cdots+\frac{|u_{\theta_{n-1}}|^2}{(\sin
r\sin\theta_1\cdots\sin\theta_{n-2})^2},
\end{align*}
and
\begin{align*}
D^2\rho(r)&=\rho^{''}(r)dr\otimes dr+\rho'(r)\sin r\cos
rd\theta_1\otimes d
\theta_1\\
&\quad+\rho'(r)\sin r\cos
r\sin^2\theta_1d\theta_2\otimes d\theta_2\\
&\quad+\cdots+\rho'(r)\sin r\cos
r\sin^2\theta_1\cdots\sin^2\theta_{n-2}d\theta_{n-1}\otimes
d\theta_{n-1},
\end{align*}
i.e.,
\begin{align*}
&\quad D^2\rho(\nabla u,\nabla\bar{u})&\\
&=\rho^{''}(r)|u_r|^2+\rho'(r)\frac{\cos
r}{\sin^3r}(|u_{\theta_1}|^2+
\cdots+\frac{|u_{\theta_{n-1}}|^2}{(\sin\theta_1\cdots\sin\theta_{n-2})^2})\\
&=(1-(n-1)\frac{\cos r}{\sin^nr}\int_0^r\sin^{n-1}sds)|u_r|^2\\
&\quad+\frac{\cos
r}{\sin^{n}r}\int_0^r\sin^{n-1}sds(\frac{|u_{\theta_1}|^2}{\sin^2r}+\cdots+\frac{|u_{\theta_{n-1}}|^2}{(\sin
r\sin\theta_1\cdots\sin\theta_{n-2})^2})\\
&\leq g(\nabla u,\nabla\bar{u}).
\end{align*}
Following the analysis in section 3 and 4, we can easily obtain
Theorem 11 and 12.

\begin{Thm}
Consider the Schr$\ddot{o}$dinger equations (3) on $M=S^n_{+}$.
Assume there exists a constant $\kappa\geq3$ such that
$$
sF(s)\geq \kappa G(s),\ \ \forall s\geq0.
$$
Then any solution satisfying (4) with $E_0<0$ blows up in finite
time.
\end{Thm}

\begin{Thm}
Consider the Schr$\ddot{o}$dinger equation
\begin{align*}
\left\{\begin{array}{ll}
iu_t=\Delta u+|u|^{p-1}u,\ \ \textrm{on} \ S^n,\\
u(0)=u_0\in H^1(S^n).
\end{array}
\right.
\end{align*}
For $p\geq5$, if $u_0(x)=u_0(r)$, where
$r=\textrm{dist}(\textrm{N},x)$ and $u_0(r)$ is an asymmetric
function at $r=\frac{\pi}{2}$ with $u(\frac{\pi}{2})=0$ and $E_0<0$,
then the asymmetric solution satisfying $(4)$ blows up in finite
time.
\end{Thm}

We next compute the case $M=H^n(-1)$. As above, for $0<r<\infty$,
$0\leq\theta_1,\cdots,\theta_{n-1}<2\pi$, we have
\begin{align*}
g=dr^2+\sinh^2r(d\theta_1^2+\sin^2\theta_1d\theta_2^2+\cdots+\sin^2\theta_1\cdots\sin^2\theta_{n-2}d\theta_{n-1}^2),
\end{align*}
and
\begin{align*}
\Delta \rho(r)=\frac{(\rho'(r)\sinh^{n-1}r)'}{\sinh^{n-1}r}.
\end{align*}
Solving the ODE
$$\Delta \rho(r)=1\ \ on \ \ (0,\pi/2]$$ we get the desired
positive solution
\begin{align*}
\rho(r)=\int_0^r(\int_0^s\sinh^{n-1}\tau
d\tau)(\sinh^{n-1}s)^{-1}ds.
\end{align*}
By defining $\rho(0)=0$ we see that $\rho\in C^4(S^n_{+})$.

From (12) and (14) we get for all $u\in C^1(S^n_{+})$
\begin{align}
g(\nabla
u,\nabla\bar{u})&=|u_r|^2+\frac{|u_{\theta_1}|^2}{\sinh^2r}+\cdots+\frac{|u_{\theta_{n-1}}|^2}{(\sinh
r\sin\theta_1\cdots\sin\theta_{n-2})^2}
\end{align}
and
\begin{align*}
D^2\rho(r)&=\rho^{''}(r)dr\otimes dr+\rho'(r)\sinh r\cosh
rd\theta_1\otimes d
\theta_1\\
&\quad+\rho'(r)\sinh r\cosh
r\sin^2\theta_1d\theta_2\otimes d\theta_2\\
&\quad+\cdots+\rho'(r)\sinh r\cosh
r\sin^2\theta_1\cdots\sin^2\theta_{n-2}d\theta_{n-1}\otimes
d\theta_{n-1},
\end{align*}
i.e.,
\begin{align}
&\quad D^2\rho(\nabla u,\nabla\bar{u})&\\
&=\rho^{''}(r)|u_r|^2+\rho'(r)\frac{\cosh
r}{\sinh^3r}(|u_{\theta_1}|^2+\frac{|u_{\theta_2}|^2}{\sin^2\theta_1}+
\cdots+\frac{|u_{\theta_{n-1}}|^2}{(\sin\theta_1\cdots\sin\theta_{n-2})^2})\nonumber\\
&=(1-(n-1)\frac{\cosh r}{\sinh^nr}\int_0^r\sinh^{n-1}sds)|u_r|^2\nonumber\\
&\quad+\frac{\cosh
r}{\sinh^{n}r}\int_0^r\sinh^{n-1}sds(\frac{|u_{\theta_1}|^2}{\sinh^2r}+\cdots+\frac{|u_{\theta_{n-1}}|^2}{(\sinh
r\sin\theta_1\cdots\sin\theta_{n-2})^2}).\nonumber
\end{align}
\textbf{Claim}:
$$
\frac{1}{n}<\frac{\cosh
r}{\sinh^{n}r}\int_0^r\sinh^{n-1}sds<\frac{1}{n-1}, \ \ \forall r>0.
$$
\begin{proof}
For $r>0$, we have
$$
\frac{\cosh
r}{\sinh^{n}r}\int_0^r\sinh^{n-1}sds>\frac{1}{\sinh^{n}r}\int_0^r\sinh^{n-1}s\cosh
sds=\frac{1}{n},
$$
and the left-hand side is proved. For the right-hand side, we write
$$
\phi(r)=\frac{\cosh r}{\sinh^{n}r}\int_0^r\sinh^{n-1}sds.
$$
Assume that $\phi(r)$ achieves its maximum at $r=r_0$,  then we have
\begin{align*}
\phi'(r_0)&=\sinh^{n-1}r_0(\cosh
r_0\sinh^nr_0+\sinh^2r_0\int_0^{r_0}\sinh^{n-1}sds\\
&\quad-n\cosh^2r_0\int_0^{r_0}\sinh^{n-1}sds)/\sinh^{2n}r_0\\
&=0,
\end{align*}
which gives that
$$
\int_0^{r_0}\sinh^{n-1}sds=\frac{\cosh
r_0\sinh^nr_0}{n\cosh^2r_0-\sinh^2r_0}.
$$
Therefore at $r=r_0$, we have
\begin{align*}
\phi(r_0)&=\frac{\cosh^2r_0}{n\cosh^2r_0-\sinh^2r_0}\\
&=\frac{\cosh^2r_0}{(n-1)\cosh^2r_0+1}<\frac{1}{n-1}.
\end{align*}
The upper bound $\frac{1}{n-1}$ is sharp due to the obvious fact
$$
\lim_{r\rightarrow\infty}\phi(r)=\frac{1}{n-1}.
$$
\end{proof}

Comparing (18) and (19) with the help of the \textbf{Claim}, we
obtain that
$$
D^2\rho(\nabla u,\nabla\bar{u})\leq\frac{1}{n-1}g(\nabla
u,\nabla\bar{u}), \ \ \forall u\in C^2(H^n(-1)).
$$
Following the analysis in section 3 and 4, we can easily obtain
Theorem 13 below on the hyperbolic space $H^n(-1)$.

\begin{Thm}
Consider the Schr$\ddot{o}$dinger equations (2) on $M=H^n(-1)$.
Assume there exists a constant $\kappa\geq 1+\frac{2}{n-1}$ such
that
$$
sF(s)\geq \kappa G(s),\ \ \forall s\geq0.
$$
Then any solution satisfying (4) with $E_0<0$ blows up in finite
time.
\end{Thm}

If $F(s)=s^\frac{p-1}{2}$ where $p>1$, the condition $sF(s)\geq
\kappa G(s)$ on $[0,+\infty)$ for some $\kappa\geq 1+\frac{2}{n-1}$
is equivalent to $p\geq 1+\frac{4}{n-1}$. We state this result as
Theorem 14.

\begin{Thm}
Consider the Schr$\ddot{o}$dinger equation
\begin{align*}
\left\{\begin{array}{ll}
iu_t=\Delta u+|u|^{p-1}u,\ \ \textrm{on}\ H^n(-1),\\
u(0)=u_0\in H^{1},
\end{array}
\right.
\end{align*}
where $p\geq 1+\frac{4}{n-1}$. Then any solution satisfying (4) with
$E_0<0$ blows up in finite time.
\end{Thm}

We remark that the similar result is also true for the complete
warped product manifold $M:=\mathbb{R}_+\times \mathbb{B}^{n-1}$
with the metric $g=dr^2+k(r)^2ds^2$. Here $(\mathbb{B}^{n-1}, ds^2)$
is a closed manifold of dimension $n-1$ and $k(r)$ is a non-negative
smooth function in $[0,\infty)$  with $k(0)=0$, $k(r)>0$ for all
$r>0$, and $k'(0)>0$. Assume that $(M,g)$ has bounded geometry. Then
we choose the function
$$
\rho=\rho(r):=\int_0^rk^{-(n-1)}(s)\left(\int_0^s
k(\tau)^{n-1}d\tau\right)ds
$$
as the weight function in
$$
J(u)=\int_M \rho |u|^2.
$$
As showed in page 31 of \cite{Ch}, we have
$$
\Delta \rho=1, \ \ \rho(0)=0, \ \ \rho'\geq 0,
$$
with uniform bounded Hessian
$$
D^2\rho\leq cg
$$
for some positive constant $c$. Then, using the similar argument, we
have the following result:

\begin{Thm}
Consider the Schr$\ddot{o}$dinger equations (2) on the warped
product $M$ as above. Assume there exists a constant $\kappa\geq
2c+1$ such that
$$
sF(s)\geq \kappa G(s),\ \ \forall s\geq0.
$$
Then any solution satisfying (4) with $E_0<0$ blows up in finite
time.
\end{Thm}

 Finally, we remark that our  blow-up result can be extended to the
 following Schr$\ddot{o}$dinger equations with harmonic potential
 on the warped product
 manifold $M$ above:
\begin{align*}
\left\{\begin{array}{ll}
iu_t=\Delta u+\rho u+|u|^{p-1}u,\ \ \textrm{on} \ M,\\
u(0)=u_0\in H^1(M).
\end{array}
\right.
\end{align*}
However, we prefer not to give the detailed statements here. Using
the similar method, we can also set up the blow-up result for
Klein-Gordon equations on Riemannian manifolds. We believe our
method can also be used to some other evolution systems (see
\cite{Bo4}).

\end{document}